\documentclass[12pt,twoside]{amsart}
\usepackage{amssymb}
\usepackage{amscd}

\numberwithin{equation}{section} 

\newtheorem{theorem}{Theorem}[section]
\newtheorem{lemma}[theorem]{Lemma}
\newtheorem{proposition}[theorem]{Proposition}
\newtheorem{corollary}[theorem]{Corollary}

\theoremstyle{definition}
\newtheorem{definition}[theorem]{Definition}
\newtheorem{remark}[theorem]{Remark}

\newcommand\Spec{\operatorname{Spec}}

\newcommand\Hom{\operatorname{Hom}}
\newcommand\Ext{\operatorname{Ext}}

\newcommand\length{\operatorname{length}}

\newcommand{\Proj}{\operatorname{Proj}}

\newcommand{\cI}{{\mathcal I}}

\newcommand{\cF}{{\mathcal F}}

\newcommand{\cO}{{\mathcal O}}
\newcommand{\cU}{{\mathcal U}}

\newcommand {\ZZ}{\mathbb{Z}}

\newcommand {\PP}{\mathbb{P}}
\newcommand {\AAA}{\mathbb{A}}

\begin{document}
\title[Families of space curves with large cohomology]{Families of space
curves with
large cohomology}

\author[Nadia Chiarli, Silvio Greco, Uwe Nagel]{Nadia Chiarli$^*$, Silvio
Greco$^*$,
Uwe Nagel}

\address{Dipartimento di Matematica, Politecnico di Torino,
I-10129 Torino, Italy}
\email{nadia.chiarli@polito.it}
\address{Dipartimento di Matematica, Politecnico di Torino, I-10129 Torino,
Italy}
\email{silvio.greco@polito.it}
\address{Department of Mathematics,
University of Kentucky, 715 Patterson Office Tower, Lexington, KY
40506-0027, USA}
\email{uwenagel@ms.uky.edu}




\thanks{$^*$ Supported by MIUR}


\begin{abstract} We investigate space curves with large cohomology. To this
end we introduce curves of subextremal type. This class includes all
subextremal curves.
Based on geometric and numerical characterizations of curves of subextremal
type, we show
that, if the cohomology is ``not too small,'' then they can be
parameterized by the union
of two generically smooth irreducible families; one of them corresponds to the
subextremal curves. For curves of negative genus, the general curve of each
of these
families is also a smooth point of the support of an irreducible component
of the
Hilbert scheme. The two components have the same (large) dimension and meet
in a
subscheme of codimension one.
\end{abstract}


\maketitle

\tableofcontents


\section{Introduction} \label{section-intro}

In this note we study space curves of degree $d$ and (arithmetic) genus $g$
that have
large cohomology. Since in this case the cohomology puts only little
restrictions on the
curves to deform, one expects that such curves form large families. In fact,
Martin-Deschamps and Perrin have shown that among the curves $C$ with fixed
$d$ and $g$,
there are curves that maximize the Rao function $h^1 (\cI_C (j))$  for all
$j \in \ZZ$.
Such curves are called extremal curves. They also showed in
\cite{MDP3} that the extremal curves form a family whose closure in the
Hilbert scheme
$H_{d, g}$ of locally Cohen-Macaulay curves is, topologically, a
generically smooth
component of
$H_{d, g}$.

If one excludes the extremal curves, Nollet \cite{N} showed that among
the remaining
curves, there are again curves that maximize $h^1 (\cI_C (j))$  for
all $j \in \ZZ$.
These curves are called subextremal. However, one cannot continue in
this fashion. Among
the curves that are neither extremal nor subextremal, there is no
curve that maximizes
the Rao functions in  {\em all}  degrees. This motivates our
definition of curves of
subextremal type. These are curves that have the same Rao function as
the subextremal
curves in all degrees $j = 1,\ldots,d-3$. Each such curve is contained
in a unique
quadric that is either reducible or not reduced.

It turns out that the curves of subextremal type can be parameterized by
two irreducible
and generically smooth families that have the same dimension. One of them
corresponds to
subextremal curves. The curves in the other family have the property that
each of them is
contained in a quadric that is not reduced, i.e.\ in a double plane.
Furthermore, we show
that if $g < 0$, then the closure of each of the two families in $H_{d, g}$ is,
topologically, a generically smooth irreducible component of $H_{d, g}$.
The two
components meet in a subscheme of codimension one that corresponds to the
subextremal
curves that are contained in  a double plane.

The paper is organized as follows:

Section 2 contains some preliminary results. After recalling the
definitions and
characterizations of subextremal and extremal curves, we establish some
useful tools. We
discuss the residual sequence of a curve $C$ with respect to a hyperplane
$H$ that
contains a planar subcurve $C' \subset C$. This sequence determines a
zero-dimensional
subscheme $Z \subset H$. For curves of subextremal type, $Z$ turns out to
be contained in
a conic. This puts heavy restrictions on $Z$ which are pointed out at the
end of Section \ref{prel}.

The following section is devoted to the structure of curves of subextremal
type (Theorem
\ref{characterize_subextremal_type}). In particular, we show that a curve
of subextremal
type can be characterized by the values of its Hilbert function in degree
two and three.
Geometrically, it is distinguished by the fact that it contains a planar
subcurve of
degree $d-2$ and that the residual curve $C'$ is a planar conic. This is
used to
completely  describe the Rao functions of curves of subextremal type (Theorem
\ref{rao_functions_subextremal_type}). In addition to the degree $d$ and
genus $g$, each
such Rao function is determined by an integer $b$ where $b$ can take only
finitely many
values.

In Section \ref{HRset}, we investigate numerical invariants of a curve $C$
of subextremal
type. Whereas the postulation character of $C$ depends only on its degree
and genus, its
graded Betti numbers are determined by the triple $(d, g, b)$ and depend
also on $b$.
This is a consequence of results in \cite{CGN6} which we also apply to
determine the
defining equations of the curves that are not subextremal.

In Section \ref{Family_SET} we begin our study of families of curves of
subextremal type.
More precisely, we exhibit two families which are distinct if the
cohomology is not too
small. The first one parameterizes the curves of subextremal type that are
contained in a
double plane. This family can be stratified according to the Rao function
of its curves.
The curves with minimal Rao function form an open and generically smooth
subfamily. The
second family is formed by the subextremal curves. Both families give rise
to irreducible
and generically smooth subschemes of the Hilbert scheme $H_{d, g}$. They
have the same
(large) dimension, which is explicitly computed.

In Section \ref{sec-hilb} we give a geometric description of the general
curve in each of
the two families. We use it to show that the closure of the two subschemes
of $H_{d, g}$
corresponding to the two families of curves of subextremal type are
actually the support
of irreducible  components of $H_{d, g}$, provided $g < 0$. Thus, in this
case the
Hilbert scheme $H_{d, g}$ contains besides the component that parameterizes
extremal
curves two further components, one is smooth at the general subextremal
curve, the other
is smooth at the general curve of subextremal type that is contained in a
double plane.


\section{Preliminary results and background} \label{prel}

We collect here some results that will be used later on.
\medskip


\noindent
{\bf 2.1. Standing Notation}

\begin{itemize}

\item $K$:  algebraically closed field of characteristic zero.

\item $\PP^n$ the $n$-dimensional projective space over $K$.

\item For a closed subscheme $X \subseteq \PP^n$, $h_X$ denotes the
Hilbert function of $X$ and $\partial h_X$ denotes the first difference of
$h_X$, i.e.\
$\partial h_X (j) = h_X (j) - h_X (j-1)$.

\item If $X \subseteq \PP^n$ is a closed subscheme, then
$\mathcal I_X \subseteq \mathcal O_{\PP^n}$
denotes the ideal sheaf of $X$ and $I_X \subseteq K[X_0,\ldots,X_n]$
denotes the (saturated)
homogeneous ideal of $X$.

\item We agree that the empty subscheme of $\PP^n$ has degree $0$.

\item $C \subseteq \PP^3$: non-degenerate, projective curve of degree $d$
and arithmetic genus $g$, where curve  means a pure $1$-dimensional
 projective subscheme (i.e. without $0$-dimensional components); in
particular $C$ is
locally Cohen-Macaulay.

\item $\Gamma$: general hyperplane section of $C$.

\item If $C$ is a curve,
the function $\rho_C(j):= h^1(\mathcal I_C(j))$ ($j\in \ZZ$) is called the
{\it Rao function of $C$}.
\end{itemize}

\bigskip


\noindent {\bf 2.2. Extremal and subextremal curves}
\medskip

Now we recall some results on curves having large cohomology, which were
one of the
starting points for our investigation. These curves were studied by
Martin-Deschamps and
Perrin \cite{MDP2}, Ellia \cite{E},  and Nollet \cite{N}.

Martin-Deschamps and Perrin in  \cite {MDP2} proved that for $d\geq 2$ the
Rao function
of $C$ satisfies the inequality $ \rho_C(j) \leq \rho^E(j)$, where $\rho^E:
\ZZ \to \ZZ$
is the function defined by:

\medskip
{\flushleft

\space \space $$\rho^E(j) : =
\left \{\begin{array}{*3l} 0&{\rm if}& j \leq - \binom{d-2}{2} + g\\
\\ \binom{d-2}{2} - g + j&{\rm if}& - \binom{d-2}
{2} + g \leq j \leq 0\\ \\ \binom{d-2}{2}
- g&  {\rm if}& 0 \leq j \leq {d-2}
\\ \\ \binom{d-1}{2} - g -
j&{\rm if}& d-2 \leq j \leq \binom{d-1}{2} - g\\ \\
0&{\rm if}& \binom{d-1}{2} - g \leq j.
\end{array} \right.$$}

\medskip

 A non-degenerate curve
$C\subseteq\PP^3$ such that $\rho_C(j) = \rho^E(j)$
for every $j \in \ZZ$, is called {\it extremal} (see \cite{MDP2}).

\medskip

For extremal curves, the following characterization follows by the results
in \cite{MDP3}
and \cite {E}:

\begin{theorem}\label{extremal_curve_th} Suppose $d\geq 5$. Then the
following are equivalent:
\begin{itemize}
\item[(a)] $C$ is extremal;
\item[(b)] $C$ contains a planar subcurve of degree $d-1$;
\item[(c)] $C$ is contained in two independent quadrics.

\end{itemize}
\end{theorem}

The extremal curves form an interesting family of large dimension.

\begin{theorem}\label{family extremal curves} Assume $d \ge 6$ and $g \le
\frac 12 (d-3)(d-4)+1$. Then the extremal curves of degree $d$ and genus $g$ form an
irreducible generically smooth family $\mathcal F_{EX}$ of dimension
$$
2a+4+ \frac 12(d-1)(d+2)
$$
where $a:= \binom{d-2}{2} - g$ is the maximum value of the Rao function.
\end{theorem}

For the proof see \cite {MDP3}, where also the other values of $d$ and $g$ are
considered.

\medskip

In a subsequent paper (\cite{N}),  Nollet  proved that, for $d\geq 5$, if
$C$ is not
extremal, then
 $\rho_C(j) \le \rho^{SE}(j)$, where $\rho^{SE}: \ZZ \to \ZZ$
is the function defined by:

\medskip
{\flushleft \space \space $$\rho^{SE}(j) : =
\left \{\begin{array}{*3l} 0&{\rm if}& j < g - \binom{d-3}{2} \\
\\ \binom{d-3}{2} - g + j&{\rm if}& g - \binom{d-3}
{2} + 1 \leq j \leq 0\\ \\ \binom{d-3}{2}
- g + 1 &  {\rm if}& 1 \leq j \leq {d-3}
\\ \\ \binom{d-2}{2} - g + 1 -
j&{\rm if}& d-3 \leq j \leq \binom{d-2}{2} - g\\ \\
0&{\rm if}& \binom{d-2}{2} - g + 1 \leq j.
\end{array} \right.$$}

\medskip

 A non-degenerate curve
$C\subseteq\PP^3$ such that $\rho_C(j) = \rho^{SE}(j)$ for every $j \in
\ZZ$, is called
{\it subextremal} (see \cite{N}).  Subextremal curves are classified in
\cite{N}.
We will see that the subextremal curves of degree $d$ and genus $g$ form a
smooth
irreducible  family $\mathcal F_{SE}$ of dimension
$2r+6+\frac{(d-2)(d+1)}{2}$, provided $r\ge 3$, where $r:= \binom{d-3}{2}
+1 - g$ is the
maximum value of the Rao function (cf.\ Theorem \ref{families2}).

\bigskip


\noindent {\bf 2.3. Residual sequences}
\medskip

Assume that $C$ contains a planar subcurve $D$ of degree $d-\delta \le d$
spanning a plane $H$ and let $\ell \in R$
be a linear form defining $H$. Let $C'$ be the {\it residue of $C$ with
respect to $H$}, namely $\mathcal I_{C'} := \mathcal I_C: \mathcal I_H$,
and let $Z \subseteq H$ be the {\it residue of
$C\cap H$ with respect to $D$},
namely $\mathcal I_{Z,H} := \mathcal I_{C\cap H,H}: \mathcal I_{D,H}$.
Let $g'$ denote the arithmetic genus of $C'$.

\medskip

\begin{proposition}\label{residual} With the above notation we have:

\begin{itemize}
\item[(i)] $\mathcal I_{Z,H}(\delta - d)$ is isomorphic to
$\mathcal I_{C\cap H,H}$, via the multiplication by an equation of $D$;
\item[(ii)] there exists an exact sequence
(called {\it residual sequence with
respect to $H$}):
$$
0 \to \mathcal I_{C'}(-1) \to \mathcal I_C \to
\mathcal I_{Z,H}(\delta-d) \to 0
,$$
where the first map is the multiplication by $\ell$;
\item[(iii)] $C'$ is a curve of degree $\delta$;
\item[(iv)] $Z$ is either empty or zero-dimensional;
\item[(v)] $\deg (Z) = \binom {d-\delta-1}{2}-g+g'+\delta - 1$;
\item[(vi)] $Z$ is a subscheme of $C'\cap H$.
\end{itemize}
\end{proposition}
\begin{proof}

(i) - (iv) are straightforward  and (vi) follows from \cite{CGN4}, Lemma 2.8.

(v) If $Z$ is non-empty,  the residual sequence provides by considering Euler
characteristics:
$$
\begin{array}{rcl}
\deg Z& = &\chi(\mathcal O_{\PP^2}(\delta-d) - \chi(\mathcal I_Z(\delta-d))
\\ \\
&=& \binom {d-\delta-1}{2} + \chi(\mathcal I_{C'}(-1)) - \chi(\mathcal I_C)
\\ \\
&=&\binom {d-\delta-1}{2} +\chi(\mathcal O_{\PP^3}(-1)) - \chi(\mathcal
O_{C'}(-1))-\chi(\mathcal O_{\PP^3}) + \chi(\mathcal O_C)
\end{array}
$$
and the conclusion follows by a straightforward computation. If $Z$ is
empty we have
$\mathcal I_{Z,H} = \mathcal O_H$ and the conclusion follows by a similar
argument.
\end{proof}

\medskip


\noindent {\bf 2.4. Zero-dimensional subschemes of a conic}
\medskip

Let $E \subseteq \PP^2 = \Proj\, (S)$ be a conic where $S = K[y, z, t]$.
Let $W\subseteq
E$ be a zero-dimensional closed subscheme of degree $r$. Then it is easy to
see that $W$
satisfies:
\begin{itemize}

\item There is an integer $b$, with
$0\le b \le \frac{r-1}{2}$, such that $\partial h_W = h_b $, where $h_b:
\ZZ \to \ZZ$ is
the function:
\medskip
{\flushleft \space  $$h_b (j) : =
\left \{\begin{array}{*3l} 0&{\rm if}& j <0\\
1&{\rm if}&j = 0\\
2&{\rm if}&1 \le j \le b \\
1&  {\rm if}& b+1 \leq j \le r-1-b\\
0&{\rm if}& j > r-1-b.
\end{array} \right.$$}

\medskip
\item If $b < \frac {r-1}{2}$, then there exists a closed subscheme
$W'\subseteq W$,
which is collinear, of degree $r-b$.
\medskip

\item If $b = \frac {r}{2}-1$, then either $W$ is a complete intersection
or again
there exists a closed subscheme $W'\subseteq W$, which is
collinear, of degree $r-b$.
\medskip

\item $W$ is collinear if and only if
 $b=0$.

\medskip

\item
 $I_W$ can have at most three minimal generators and
there are the following possibilities:

\underline {Case 1}. $W$ is collinear. Then $I_W$ is a complete
intersection of type
$(1,r)$ and its minimal free resolution has the form:

$$0\to S(-1-r)\to S(-1)\oplus S(-r)\to I_W\to 0.$$

\medskip

\underline {Case 2}. $W$ is a complete intersection of type $(2,\frac
{r}{2})$. Then its
minimal free resolution has the form:

$$0\to S(-2-\frac {r}{2})\to S(-2)\oplus S(-\frac {r}{2})\to I_W\to 0.$$

\medskip
\underline {Case 3}. $W$ is not a complete intersection. Then $I_W$ has
exactly three
minimal generators of degree $2, b+1, a+1$, with $2\le b+1\le a+1$ and  $a
= r-b-1$,
where $b$ the integer which defines the Hilbert function of $W$ (see
above). Moreover,
the minimal free resolution of $I_W$ has the form:

\medskip

{\flushleft
\space \space $$
 \begin{array}{*9c} &&&&S(-2)&&&&\\
&&S(-b-2)&&\oplus &&&&\\
0&\to &\oplus &\overset \varphi \to &S(-b-1)& \to &I_W&\to &0. \\
&&S(-a-2)&&\oplus &&&&\\
&&&&S(-a-1)&&&&
\end{array} $$}
The degree matrix of the Hilbert-Burch matrix representing $\varphi$ is
$$\left [\begin{array}{*2c} b&a\\
1&a-b+1\\
b-a+1&1
\end{array}\right ].$$

Furthermore, the ideal of $W$ is $I_W=I(\varphi)$, where $I(\varphi)$
denotes the ideal
generated by the maximal minors of the matrix.

\end{itemize}


\section{Structure theorem for curves of subextremal type}
\label{Structure}

In this section we consider a class of curves with large cohomology. It
turns out that
they have a rather particular structure. We use it to determine all
occurring Rao
functions among these curves.

\begin{definition}\label{subextremal_type} A non-degenerate curve
$C\subseteq\PP^3$ of degree $d$ and genus $g$ is said to be {\it of
subextremal type} if
$d\ge 5$ and $\rho_C(j) = \binom{d-3}{2}-g+1$ for $1 \le j \le d-3$.
\end{definition}

\medskip

Note that, by the results of \S \ref{prel}, a subextremal curve is of
subextremal type
and that a curve of subextremal type is not extremal.

From now on, let $r := \binom{d-3}{2}-g+1$.

\bigskip

For curves of subextremal type, there is the following structure theorem:

\medskip
\begin{theorem} \label{characterize_subextremal_type}
Let $C\subseteq\PP^3$ be a non-degenerate curve of degree $d\ge 7$. Then
the following are equivalent:
\begin{itemize}
\item[(i)] $C$ is of subextremal type;
\item[(ii)] $h^0(\mathcal I_C(2))=1$ and $h^0(\mathcal I_C(3))=5$ (that is,
$I_C$
has one minimal generator in  degree $2$ and one in degree $3$);
\item[(iii)] $C$ is contained in a unique quadric and
$\partial h_\Gamma: 1\; 2\; 2\; 1\; \dots\; 1\; 0\; \rightarrow$;
\item[(iv)] $C$ contains a planar subcurve of degree $d-2$ and
the residual curve $C'$ is a planar curve of degree $2$.

\end{itemize}
\end{theorem}

\begin{proof} (i) $\Rightarrow$ (ii). Since $C$ is not extremal we have
$\partial h_\Gamma(2)
\geq 2$. It follows $h^1(\mathcal I_\Gamma (2))\leq d-5$, whence
$h^1(\mathcal I_\Gamma
(j))\leq d-3-j$ for $2\leq j\leq d-5$ and $h^1(\mathcal I_\Gamma (j))=0$
for $j>d-5$.
Hence, with an argument as in \cite {CGN1} (proof of Theorem 2.1, step 2),
we get, for
$2\leq j \leq d-5$,

$$h^2(\mathcal I_C(j))\leq \sum_{t \ge j+1} h^1(\mathcal I_\Gamma (t))\leq
\binom {d-3-j}{2}.
$$
In particular $h^2(\mathcal I_C(j))=0$ for $j\geq d-5$.

Moreover by Riemann-Roch we have, for $1\leq j\leq d-3$,
$$\begin{array} {rcl}
h^0(\mathcal I_C(j))&=&h^0(\mathcal O_{\PP^3}(j))- h^0(\mathcal
O_C(j))+h^1(\mathcal I_C(j))\\
&=&\binom {j+3}{3}-[dj-g+1+h^2(\mathcal I_C(j))]+r\\
&=& \binom {j+3}{3}-dj- h^2(\mathcal I_C(j))+\binom {d-3}{2}.
\end{array}
$$

It follows that $h^0(\mathcal I_C(2))\geq 1$ and $h^0(\mathcal I_C(3))\geq
5$. Since $C$
is not extremal, we obtain $h^0(\mathcal I_C(2))=1$. Moreover, if $h^0(\mathcal
I_C(3))>5$, then the exact sequence:
$$
0\to H^0(\mathcal I_C(2))\to H^0(\mathcal I_C(3))\to H^0(\mathcal
I_\Gamma(3))\to \cdots
$$
provides $h^0(\mathcal I_\Gamma(3))>4$, that is $h_\Gamma(3)\leq 5$. Since
$h_\Gamma(2)=
5$, this implies $5=\deg \Gamma =d$, a contradiction.

(ii) $\Rightarrow$ (iii). From the exact sequence
$$
0\to H^0(\mathcal I_C(1))\to H^0(\mathcal I_C(2))\to H^0(\mathcal
I_\Gamma(2))\to \cdots
$$
we have $h^0(\mathcal I_\Gamma(2))\geq 1$ and from the exact sequence:

$$ 0\to H^0(\mathcal I_C(2))\to H^0(\mathcal I_C(3))\to H^0(\mathcal
I_\Gamma(3))\to \cdots$$
we have $h^0(\mathcal I_\Gamma(3))\geq 4$.

If $h^0(\mathcal I_\Gamma(2))\geq 2$, then  $C$ is extremal, a contradiction.
Then we must have $h^0(\mathcal I_\Gamma(2))=1$ whence
$\partial h_\Gamma(1)= \partial h_\Gamma(2) = 2$.

Moreover we have $\partial h_\Gamma(3)=1$ (as in the proof of (i)
$\Rightarrow$ (ii)) and the conclusion follows.

(iii) $\Rightarrow$ (iv). Since $d\ge 7$, by \cite {CGN3}, Corollary 4.4,
$C$ contains a subcurve of degree $d-2$ spanning a plane $H$.
The residual sequence with respect to $H$
is:
$$
0 \to \mathcal I_{C'}(-1) \to \mathcal I_C \to \mathcal I_{C\cap H,H} \to 0
,$$
where $C'$ is a curve of degree $2$ and the first map is the
multiplication by a linear form defining $H$ (see Proposition \ref{residual}).
The sequence above provides the exact sequence:
$$
0 \to H^0(\mathcal I_{C'}(1)) \to H^0(\mathcal I_C(2))
\to H^0(\mathcal I_{C\cap H,H}(2)).
$$
Since $d-2 > 2$, we get $h^0(\mathcal I_{C\cap H,H}(2)) = 0$, whence
$h^0(\mathcal I_{C'}(1)) = h^0(\mathcal I_C(2)) = 1$;
thus $C'$ is a planar curve.

(iv) $\Rightarrow$ (i). Let $D$ be the planar subcurve of degree $d-2$. Let
$H$ be the
plane that is spanned by $D$. Let $Z \subseteq H$ be the residual scheme of
$C \cap H$
with respect to $D$. Then by Proposition \ref{residual}(iv) (with $\delta =
2$ and $g' =
0$) we have $\deg Z = r$ and the residual sequence with respect to $H$ becomes:
$$
0 \to \mathcal I_{C'}(-1) \to \mathcal I_C \to \mathcal I_{Z,H}(2-d) \to 0
.$$

Since $C'$ is a planar curve of degree $2$ we have $h^1(\mathcal I_{C'}(t))=0$
for all $t \in \ZZ$ and $h^2(\mathcal I_{C'}(t))=0$ for $t \ge 0$. Then,
for $j \ge 1$ we get:
$$
h^1(\mathcal I_C(j)) = h^1(\mathcal I_{Z,H}(2-d+j))
.$$
If $Z\ne \emptyset$, then $h^1(\mathcal I_{Z,H}(t)) = \deg (Z) = r$ for $t
\le -1$, whence
$h^1(\mathcal I_C(j)) = r$ for $1 \le j \le d-3$.

If $Z =\emptyset$, then $h^1(\mathcal I_{Z,H}(t))$ = $h^1(\mathcal
O_H(t)) = 0$ for every $t$, thus $C$  is arithmetically
Cohen-Macaulay.
\end{proof}

\begin{remark}\label{counterexamples} Theorem
\ref{characterize_subextremal_type}
is false without the assumption $d \ge 7$. Indeed consider the following
examples:
\begin{itemize}
\item[(i)] Let $C$ be a curve of type (1,4) on a smooth quadric $Q$. Then
$d = 5$ and $g = 0$, whence $r = 2$. A straightforward calculation shows that
 $\rho_C(1) = \rho_C(2) = 2$, which implies that $C$ is of subextremal type.
On the other hand it is easy to see that (ii) of Theorem
\ref{characterize_subextremal_type} does not hold for this curve.
\item[(ii)] Let $C$ be a curve of type (1,5) on a smooth quadric $Q$.
Then it is easy to see that $\partial h_\Gamma : 1\;2\;2\;1\;0
\rightarrow$, whence (iii)
of Theorem \ref{characterize_subextremal_type} holds. But it can be shown
by direct
calculations that $\rho_C(1) = 3$ and $\rho_C(2) = 4$, whence $C$ is not of
subextremal
type.
\end{itemize}
\end{remark}

The next Corollary summarizes properties of curves of subextremal type
which follow from
Theorem \ref{characterize_subextremal_type} and its proof. We state them
here for later
use.

\begin{corollary}\label{properties_subextremal_type}
Let $C$ be a curve of subextremal type of degree $d\ge 7$.
Then we have:

\begin{itemize}
\item[(i)] $C$ contains a planar subcurve $D$ of degree $d-2$
spanning a plane $H$;

\item[(ii)] the residual exact sequence with respect  to $H$ is
 \begin{equation}\label{residual_sequence}
0 \to \mathcal I_{C'}(-1) \to \mathcal I_C \to \mathcal I_{Z,H}(2-d) \to 0
\end{equation}
where $C'$ a curve of degree $2$ spanning a plane $H'$
and $Z \subseteq H$ is a closed $0$-dimensional subscheme with $\deg Z = r$;

\item[(iii)] $h^1(\mathcal I_C(j)) = h^1(\mathcal I_{Z,H}(d-2+j))$ for
$j\ge 0$;

\item[(iv)] $C$ is contained in a unique quadric $Q$ which is either the
union of
$H$ and $H'$, if $H \ne H'$, or it is the double plane $2H$, if $H = H'$.

\item[(v)] $Z \subseteq C' \cap H$. Hence if $Z \ne \emptyset$ we have
either $h^0(\mathcal I_{Z,H}(1)) \ne 0$ or $h^0(\mathcal I_{Z,H}(1)) = 0$ and
$h^0(\mathcal I_{Z,H}(2)) \ne 0$.

\item[(vi)] $Z=\emptyset$ if and only if $C$ is arithmetically Cohen-Macaulay.
\end{itemize}

\end{corollary}

\medskip

Now we completely describe the Rao functions of curves of subextremal type.


\begin{theorem}\label{rao_functions_subextremal_type}
Let $C$ be  a curve of subextremal type of degree $d \ge 7$ that is not
arithmetically
Cohen-Macaulay. Then we have:

\begin{itemize}
\item[(i)] the Rao function of $C$ is symmetric, namely:
$$
\rho_C(j) = \rho_C(d-2-j) \quad {\rm for\; all} \; j \in \ZZ;
$$
\item[(ii)] there is an integer $b$, $0\le b \le \frac{r-1}{2}$,
such that $\rho_C = \rho_b$, where $\rho_b: \ZZ \to \ZZ$ is the function
defined by:

{\flushleft \space \space $$\rho_b(j) = \left \{\begin{array}{*3l}
\rho_b(d-2-j) &{\rm if}& j \leq 0\\ \\
r &{\rm if}& 1\leq j \leq d-3\\ \\
r-1&{\rm if}& j = d-2\\ \\
r-1-2(j-d+2)&{\rm if}&d-2 \le j \le d-2+b \\ \\
r-1-b-j+d-2&  {\rm if}& d-2+b \leq j \leq d+r-3-b \\ \\
0&{\rm if}& j\ge d+r-2-b;
\end{array} \right.$$}
\item[(iii)] let $Z$ be the $0$-dimensional subscheme defined in Corollary \ref
{properties_subextremal_type}. Then the following conditions are equivalent:
\begin{itemize}
\item[(a)] $C$ is subextremal;
\item[(b)] $Z$ is collinear;
\item[(c)] $b = 0$;
\item[(d)] $\rho_C(d+r-3) > 0$;
\end{itemize}
\item[(iv)]  if the unique quadric containing $C$ is reduced, then
$C$ is subextremal.
\item[(v)] $C$ is minimal in its biliaison class if and only if $C$ is not
subextremal.

\end{itemize}
\end{theorem}

\begin{proof}
(i) Let $Q$ be the unique quadric containing $C$ (see Corollary
\ref{properties_subextremal_type}). If $Q$ is a double plane the
symmetry follows from \cite{HS}, Corollary 6.2.

If $Q$ is reduced, then by Corollary \ref{properties_subextremal_type} we
have $Q = H
\cup H'$, where $H$ is the plane containing the planar subcurve of degree
$d-2$ of $C$
and $H' \not = H$ is the plane of $C'$. Then $Z$ is contained in the line
$H\cap H'$ by
Corollary \ref{properties_subextremal_type} and hence the homogeneous ideal
$I_{Z,H}$ has
a minimal generator of degree $1$. From the residual exact sequence with
respect to $H$
it is not difficult to see  that $C$ is contained in a surface $F$ of
degree $d-1$ with
no common components with $Q$. Let $E$ be the curve linked to $C$ by the
complete
intersection $Q\cap F$. By liaison (see e.g. \cite M)) one has: $d':= \deg
E = d-2 \ge
5$, $g': = p_a(E) = g-(d-3)$, $\rho_E(j)) = \rho_C(d-3-j))$ for all $j\in
\ZZ$. It
follows that $\rho_E(0) = r > \frac 12(d'-3)(d'-4)+1-g'$, whence $E$ is
extremal by
Nollet's bound (\cite {N}). The conclusion follows by the definition of
extremal curves.

(ii) By (i) we may assume $j\ge 0$. For such values of $j$ we have, by
Corollary
\ref{properties_subextremal_type}(iii), $h^1(\mathcal I_C(j)) =
h^1(\mathcal I_Z(d-2+j))
= r - h_Z(2-d-j)$, whence $\rho_C(j) = \rho^{SE}(j)$ for $j\ge 0$ and $C$
is subextremal
by (i). The conclusion follows easily by Corollary
\ref{properties_subextremal_type}(v)
and \S \ref{prel}.

(iii) It can be proved by an argument similar to the previous one.

(iv) If $C$ is contained in a reduced quadric, then $Z$ is collinear
(see proof of (i)) and the conclusion follows from (iii).

(v) If $C$ is subextremal, then it is not minimal. Conversely assume $C$ is
not minimal;
then $C$ lies in a double plane by (iv) and $Z$ is not collinear by (iii).
Hence $C'$ is
a curve of minimal degree containing $Z$, whence $C$ is minimal by \cite
{HS}, Corollary
7.3.
\end{proof}

\begin{remark} \label{rho_b} (i) The function $\rho_b$ of Theorem
\ref{rao_functions_subextremal_type}, for
$j \ge d-2$, decreases by 2 for $b$ steps and then by 1 until it vanishes.
See the
following picture, where we put $a:=r-b-1$.

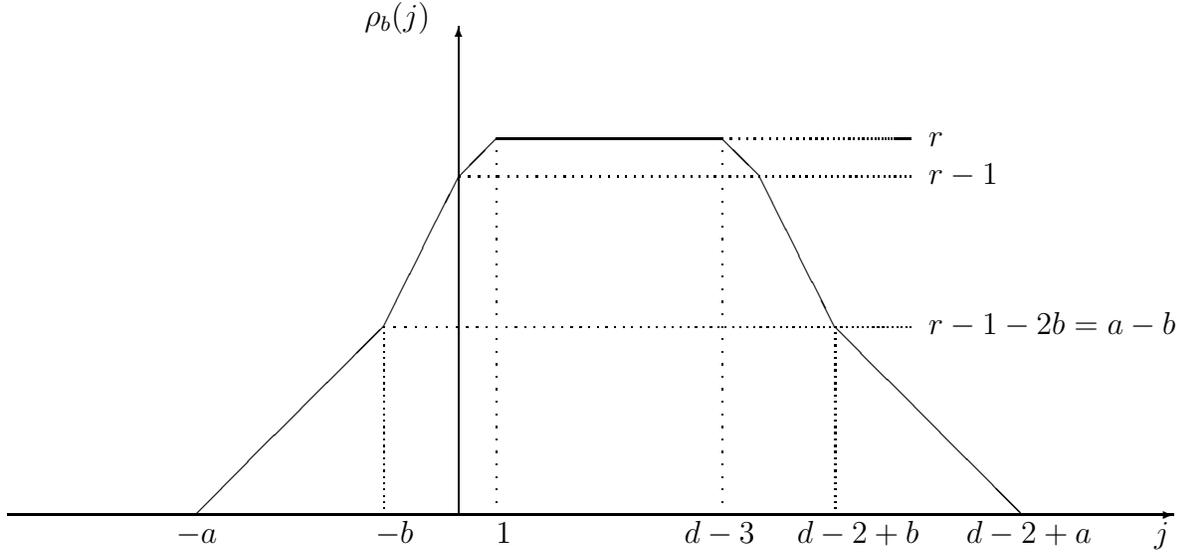
\begin{picture}(260,260)(-160,-50)
\thinlines \setlength{\unitlength}{5mm}

\put(-12,0){\vector(1,0){31}} \put(0,0){\vector(0,1){13}}



\put(1,-.75){$1$} \put(6.0,-.75){$d-3$} \put(-7.5,-.75){$-a$}
\put(-2.2,-.75){$-b$}
\put(-2.5,13){$\rho_b(j)$} \put(18.5,-.75){$j$}



\put(0,0){\qbezier[30](1,0)(1,5)(1,10)} \put(0,0){\qbezier[30](7,0)(7,5)(7,10)}
\put(0,0){\qbezier[30](10,0)(10,2.5)(10,5)}
\put(0,0){\qbezier[30](-2,0)(-2,2.5)(-2,5)}

\put(9,-.75){$d-2+b$} \put(13.5,-.75){$d-2+a$}


\put(0,0){\qbezier[40](7,10)(11.5,10)(12,10)}
\put(0,0){\qbezier[60](0,9)(8,9)(12,9)}
\put(0,0){\qbezier[70](-2,5)(8,5)(12,5)} \put(12.5,9.8){$r$}
\put(12.5,8.8){$r-1$}
\put(12.5,4.8){$r-1-2b=a-b$}
\put(-7,0){\line (1,1){5}} \put(-2,5){\line (1,2){2}}
\put(1,10){\line (1,0){6}} \put(0,9){\line (1,1){1}} \put(7,10){\line
(1,-1){1}}
\put(8,9){\line (1,-2){2}} \put(10,5){\line (1,-1){5}}
\end{picture}

Note that  $a$ and $b$ are the numbers introduced in \S 2.4 to describe the
Hilbert
function and the minimal free resolution of $I_{Z,H}$.
\medskip

(ii) We will see  that for every triple $(d, g, b)$ of integers such that
$d \geq 7,\; g
\leq \binom{d-3}{2} + 1$, and  $0\le b \le \frac{r-1}{2}$, there exists a
curve $C$ of
subextremal type having  Rao function $\rho_b$, as prescribed by Theorem
\ref{rao_functions_subextremal_type}(ii). The equations of one such curve
are described
in Theorem \ref{ideal and rao} (cf.\ Remark \ref{rem-ideals}); the
resulting families are
studied in Theorem \ref{families1}.
\end{remark}
\medskip

\section{The  ideal and numerical characters of a curve of subextremal type}
\label{HRset}

\medskip

In this section we describe information about curves of subextremal type
that we need for
studying their families. At first, we focus on curves that are not subextremal.

Let $C$ be a curve of subextremal type that is neither subextremal nor
arithmetically Cohen-Macaulay. Assume that $d \ge 7$ and $\rho_C=\rho_b$
(cf.\ Theorem
\ref {rao_functions_subextremal_type}) and set $a:= r-b-1$. Using the notation
of Corollary
\ref{properties_subextremal_type},  recall that $C$ is contained in a
double plane $2H$
and $C'\subseteq D$. We may assume $H : = \{x = 0\}$. We identify $H$ with
$\PP^2$ with
coordinates $y,z,t$ and we set
 $I_{C'} = (\phi, x)$, where $\phi \in S:=R/xR$ is a form of degree $2$ and
$I_D = (\phi h,x)$, for a suitable form $h\in S$ of degree $d-4$.

The following result provides a minimal set of generators of the
homogeneous ideal of $C$
and a minimal presentation of $M_C$. Recall that a {\em Koszul module} is a
graded
$R$-module $R/(f_1, f_2, f_3, f_4) (t)$ where $t \in \ZZ$ and $f_1, f_2,
f_3, f_4$ is a
regular sequence.

\begin{theorem} \label{ideal and rao}
With the above notation and assumptions  we have:

\begin{itemize}
\item[(i)] If $Z$ is a complete intersection (namely $I_{Z,H} = (\psi,
\phi)$), then
$$
I_C = (x^2, x\phi, \phi^2 h, \psi \phi  h +xF)
$$
where $F \in S$ is a form of degree $d-3+ \frac r2$ such that the ideal
$(\psi, \phi, F)S$ is irrelevant.

Moreover
$$M_C \cong [R/(x,\psi, \phi, F)](d-2).$$
In particular $M_C$ is a Koszul module.

\item[(ii)] If $Z$ is not a complete intersection, possibly after a
suitable choice of
coordinates and bases, we may assume
$$
A:= \left [\begin{array}{ccc}
p&y&m\\
q&n&l
\end{array}
\right]
$$
to be the transpose of a Hilbert-Burch matrix of $I_{Z,H}$ where $\phi =
{\tiny \left
|\begin{array}{cc}
y&m\\
n&l
\end{array}\right|}$,  $\deg p = b$, and $\deg q = a$. Then
$$
I_C = \left (x^2, x\phi, \phi^2 h, \phi h{\tiny \left |\begin{array}{cc}
p&m\\
q&l
\end{array}\right|} +x{\tiny \left |\begin{array}{cc}
m&F\\
l&G
\end{array}\right|}, \phi h{\tiny \left |\begin{array}{cc}
p&y\\
q&n
\end{array}\right|} + x{\tiny \left |\begin{array}{cc}
y&F\\
n&G
\end{array}\right|}\right )
$$
where $F,G \in S$ are forms such that the $2 \times {2}$ minors of the
homogeneous matrix
$$
M:= \left [\begin{array}{cccc}
p&y&m&F\\
q&n&l&G
\end{array}
\right]
$$
generate an irrelevant ideal in $S$, $\deg F =  b+d-3, \deg G = a+d-3$, and
$h \in
K[y,z,t]$ is a non-trivial form of degree $d-4$.

Moreover  $M_C$ is isomorphic to the cokernel of the map

$$
 \begin{array}{*{7}c} R(b-2)&&R(-1)&&&&R(b-1)\\
\oplus&\oplus&\oplus&\oplus&R(-d+2)&\to&\oplus\\
R(a-2)&&R(b-2)&&&&R(a-1)\\
&&\oplus&&&&\\
&&R(a-2)&&&&
\end{array} $$
\noindent defined by the matrix $\left [xE_2,M\right]$, where $E_2$ is the
$2\times2$
identity matrix.

In particular $M_C$ is minimally generated by two homogeneous elements of
degrees $1-a$
and $1-b$.
\medskip
 \end{itemize}
\end{theorem}

\begin{proof}In \cite{CGN6} it is proved that a minimal system of
generators of the
homogeneous ideal of a curve $C$ lying in  a double plane can be expressed
by using the
maximal minors of the Hilbert-Burch matrix $A$ of $Z$ and of a certain
homogeneous matrix $B$
obtained from $A$ by adding a suitable row and a suitable column.
In particular by using  the degree relations in Corollary 3.6 and Remark 4.7
of \cite{CGN6}, we may assume in case (i) that
$$B= \left [\begin{array}{ccc}
\psi&\phi&-F\\
1&0&0
\end{array}
\right]$$
where $F$ is a form of the required degree and in case (ii) that $$B= \left
[\begin{array}{cccc}
p&y&m&F\\
q&n&l&G\\
1&0&0&0
\end{array}
\right]$$ with  forms $F$ and $G$ of the required degrees.

The conclusion about the generators of $I_C$ follows immediately by the
above mentioned result.
The expressions for $M_C$ follow from \cite {CGN6}, Theorem 4.1(i).
\end{proof}
\medskip

\begin{remark} \label{rem-ideals}
(i) Note that conversely, each ideal that is defined as in the above
theorem, is
saturated and defines a curve of subextremal type. This follows from
\cite{CGN6}.

(ii) It is easy to produce equations for a specific curve with Rao function
$\rho_b$.
Given $d, g$ one computes $r = \binom{d-3}{2} - g + 1$. Let $b$ be an
integer such that
$0 \leq b \leq \frac{r}{2} - 1$. Then choose two lines $L_1 \neq L_2$ in
the plane $\{ x
= 0\}$ and, for $i = 1, 2$, subsets $Z_i \subset L_i$ of $b$ and $a = r- b
- 1$ points,
respectively,  such that the union $Z$ of $Z_1, Z_2$, and $L_1 \cap L_2$
consists of $r$
points. Let $A$ be the transpose of the Hilbert-Burch matrix of $Z$. Use
the matrix $M =
\begin{bmatrix}
A & \begin{array}{c} F \\
G
\end{array}
\end{bmatrix} $
where $F, G$ are sufficiently general forms of degree $b+d-3$ and $a + d-3$,
respectively, and, e.g.,  $h := y^{d-4}$ to create an ideal $I$ as in
Theorem \ref{ideal
and rao}. Then $I$ defines a curve $C$ of subextremal type with Rao
function $\rho_b$.
Indeed the Hilbert function of $Z$ is $h_b$ (cf. \S 2.4) and this implies
that $\rho_C =
\rho_b$ by Remark \ref{rho_b}.
\end{remark}

\begin{corollary} Let $C$ be a curve of subextremal type with $d \ge 7$.
Then
\begin{itemize}
\item[(i)] $\rho_C$ determines the structure of $M_C$, except when $r$ is
even and
$b = \frac{r}{2} -1$.
\item[(ii)] $M_C$ is a Koszul module if and only if $Z$ is a complete
intersection.
\end{itemize}
\end{corollary}

\begin{proof}
If $C$ is subextremal, this follows by \cite{N}. Otherwise, apply Theorem
\ref{ideal and
rao}.
\end{proof}

The following theorem provides the minimal free resolution of $I_C$. It is
a particular
case of \cite {CGN6}, Theorem 3.7. We write here only the modules (hence
the Betti
numbers), referring to the above mentioned result for an explicit
description of the
maps, which can be expressed in terms of the matrix $B$ given in the proof
of Theorem
\ref{ideal and rao}.
\medskip

\begin{proposition}\label{min_free_of_set} Let $C$ be as above and set $a:=
r-b-1$.
Then the minimal free resolution of $I_C$ is:

\underline {Case 1}: $Z$ is a complete intersection of type $(2,\frac
{r}{2})$ (thus $b =
\frac{r}{2} - 1$).

\medskip

$
0 \to R(-d-\frac {r}{2}-1) \to
\begin{array}{c}
R(-4)\oplus R(-d-1) \\
\oplus \\
R(-\frac {r}{2}-d+1) \\
\oplus \\
R(-\frac {r}{2}-d)
\end{array}
\to
\begin{array}{c}
R(-2)\oplus R(-3) \\
\oplus \\
R(-d)\oplus R(-\frac {r}{2}-d+2)
\end{array}
\to I_C \to 0
$

\medskip

\underline {Case 2}: $Z$ is not a complete intersection.

\medskip

{\small
 $ 0 \to
\begin{array}{c}
R(-d-b-1) \\
\oplus \\
R(-d-a-1)
\end{array}
\to
\begin{array}{c}
R(-4)\oplus R(-d-1) \\
\oplus \\
R^2(-d-b) \\
\oplus\\
 R^2(-d-a)
\end{array}
\to
\begin{array}{c}
R(-2)\oplus R(-3) \\
\oplus \\
R(-d)\oplus R(-d-b+1) \\
\oplus \\
R(-d-a+1)
\end{array}
\to I_C \to 0 $ }
\end{proposition}
\medskip

In the above result, we left out the case of subextremal curves. For these,
we have:

\begin{proposition} \label{prop-res-se}
The minimal free resolution of a subextremal curve $C$ of degree $d\ge 7$
is of the form:
\medskip

$ 0 \to R(-r-d) \to
\begin{array}{c}
R(-4)\oplus R(-d) \\
\oplus \\
R^2(-r-d+1)
\end{array}
\to
\begin{array}{c}
R(-2)\oplus R(-3) \\
\oplus \\
R(-d+1)\oplus R(-r-d+2)
\end{array}
\to I_C \to 0.
$
\end{proposition}

\begin{proof}
This follows by applying the Horseshoe Lemma to the residual sequence in
Corollary
\ref{properties_subextremal_type}(ii). Note that the resulting free
resolution is minimal
because our assumption $d \geq 7$ guarantees that no cancellation is possible.
\end{proof}


\medskip

We now turn to the computation of numerical characters. Recall the
following facts (see
\cite {MDP1}, Definition 2.3 and Proposition 2.6).

\begin{definition} The {\it postulation character} of a curve $C \subseteq
\PP^3$ is the
function $\gamma_C$ defined by
$$
\gamma_C:= -\partial ^3 h_C
$$

Observe that $\gamma_C$ and the Hilbert function $h_C$ determine each other.
\end{definition}

\begin{corollary}\label{Hilbert} Let $C$ be a curve of subextremal type of
degree $d\geq 7$.
 Then we have:

\begin{itemize}
\item[(i)]

{\flushleft
\space \space $h^2(\mathcal I_C(j)) =
\left \{\begin{array}{*3l}
\rho_C(j)-dj+g-1 &{\rm if}& j <0\\ \\
r+g-1&{\rm if}& j = 0\\ \\
\binom{d-3-j}{2}&{\rm if}&1 \le j \le d-5 \\ \\
 0&{\rm if}& j> d-5;
\end{array} \right.$}
\medskip

\item[(ii)]

{\flushleft
\space \space $h_C(j) =
\left \{\begin{array}{*3l}
dj-g+1+ \binom{d-3-j}{2}-\rho_C(j)&{\rm if}&0\leq j \leq d-5\\ \\
dj-g+1-\rho_C(j)&{\rm if}& j> d-5;
\end{array} \right.$}

\item[(iii)] the index of speciality of $C$ is $e = d-5$;

\item[(iv)] the postulation character $\gamma_C$  is given by:

{\flushleft \space \space $$\gamma_C(j)= \left \{\begin{array}{*3l}
0 &{\rm if}& j <0\\ \\
-1&{\rm if}& 0\le j \le 1 \\ \\
0&{\rm if}& j=2 \\ \\
 1&{\rm if}& j=3 \\ \\
0&{\rm if}& 4\le j\le d-1 \\ \\
\partial ^3 \rho_C &{\rm if}&  j\ge d;

\end{array} \right.$$}
\medskip

 in particular, it depends only on $\rho_C$.
\end{itemize}
\end{corollary}

\begin{proof} Adopt the notation in Corollary \ref
{properties_subextremal_type}.

(i) From the residual exact sequence with
 respect to $H$ we have the exact sequence:
$$
\dots\to H^2(\mathcal I_{C'}(j-1))\to H^2(\mathcal I_C(j))\to H^2(\mathcal
I_{Z,H}(2-d+j))
\to H^3(\mathcal I_{C'}(j-1))\to \dots.
$$

First assume $j\geq 1$. Then $H^2(\mathcal I_{C'}(j-1))=0$ and
$H^3(\mathcal I_{C'}(j-1))=H^3(\mathcal O_{\PP^3}(j-1))=0$, whence
$h^2(\mathcal I_C(j))=h^2(\mathcal I_{Z,H}(2-d+j))=
h^2(\mathcal O_{\PP^2}(2-d+j)) = h^0(\cO_{\PP^2}(d-5-j))$
and the conclusion follows in this case.

Assume now $j\leq 0$. We have $h^0(\mathcal O_C(j))=dj-g+1+h^2(\mathcal
I_C(j))$
and since
$$
h^0(\mathcal O_C(j)) =
\left \{\begin{array}{*3l}
\rho_C(j)&{\rm if}& j <0\\ \\
\rho_C(0)+1&{\rm if}& j=0
\end{array} \right.
$$
the conclusion follows, recalling that $\rho_C(0)+1=r$ by
Theorem \ref {rao_functions_subextremal_type} (ii).

(ii) We have:
$$ h_C(j)= h^0(\mathcal O_C(j))-\rho_C(j)=dj-g+1+h^2(\mathcal
I_C(j))-\rho_C(j).$$
If $j\geq 1$ the conclusion follows from (i). On the other hand we have
$h_C(0)=1$ and
the conclusion follows from the definition of $r$, since $\rho_C(0)=r-1$.

(iii) is an immediate consequences of (i) (ii) and the definitions.

(iv) Follows from (ii) and a straightforward calculation.
\end{proof}

\begin{corollary} \label{gamma}Let $C \subseteq \PP^3$ be a curve and
assume that $\gamma_C$
coincides with the postulation character of a curve of subextremal type having
the same degree and genus as $C$.

Then $C$ is of subextremal type and $\rho_C$ depends only on $\gamma_C$.
\end{corollary}

\begin{proof} By definition $\gamma_C$ determines $h_C$.
The conclusion follows from Theorem \ref{characterize_subextremal_type}
(ii) and Corollary \ref {Hilbert}(iv).
\end{proof}

\begin{remark} Recall that the {\it spectrum} of a curve $C$ is the function
$\ell_C(j) : = \partial^2 h^0(\mathcal O_C(j)) = \partial^2 h^2(\mathcal
I_C(j))$
(see \cite {Sc}) and that the {\it speciality character} of $C$ is the function
$\sigma_C := \partial \ell_C$ (see \cite{MDP1}, Definition 2.3).
Then by Corollary \ref{Hilbert}
we see that if $C$ is a curve of subextremal type then $\rho_C$ determines
$\ell_C$ and
$\sigma_C$ and conversely.
\end{remark}



\section{The family of curves of subextremal type}
\label{Family_SET}

In this section we want to show that the curves of subextremal type of
given degree
and genus form a family and to describe this family.

We consider only projective families parameterized by the closed points of
algebraic
$k$-schemes. If $X$ is a scheme and $x \in X$ we denote by $\kappa(x)$ the
residue
field of the local ring $\cO_{X,x}$.

Throughout this section we fix the integers $d$ and $g$ and we put $r:=
\binom{d-3}{2}+1
- g$. Moreover we denote by $H_{d,g}$ be the Hilbert scheme of locally
Cohen-Macaulay
curves of degree $d$ and genus $g$.
 We refer to \cite{sernesi} and \cite{MDP1} for basic information about
families and Hilbert schemes.

We begin with a result which is of course expected but still needs a proof.

\begin{lemma}\label{families_lemma1} Let $E \times \PP^3 \supseteq X \to E$
be a family of
curves of subextremal type of degree $d\ge 7$ and genus $g$ and let  $G
\times \PP^3
\supseteq Y
\to G$ be the family of quadrics of $\PP^3$. Then we have:
\begin{itemize}
\item[(a)] There is a morphism $f : E \to G$ such that  $G_{f(e)}$
is the unique quadric containing $X_e$, for each $e \in E$ (see Theorem
\ref{characterize_subextremal_type});
\item[(b)] The subset $E^{(2)} \subseteq E$ corresponding to the curves
contained in
some double plane is closed.
\end{itemize}
\end{lemma}

\begin{proof} (a) Let $U = \Spec(A) \subseteq E$ be an open affine subset
and set $B := A[x,y,z.t]$, the graded polynomial ring in $4$ variables. Then $X_U:= X
\cap (\PP^3 \times U)$ is a closed subscheme of $\PP^3 \times U = \Proj (B)$. Let $I
\subseteq B$ be the saturated homogeneous ideal of $X_U$. For each $e \in U$, $(B/I)
\otimes \kappa(e)$ is the homogeneous coordinate ring of $X_e$, whence by Theorem
\ref{characterize_subextremal_type} (iii) we have $\dim_{\kappa(e)} ((B/I)_2 \otimes
\kappa(e)) = 9$ for all $e\in U$. This implies that $(B/I)_2$ is a projective $A$-module
of rank 9, whence $I_2$ is a projective $A$-module of rank $1$. Let $J:= I_2B$. It is a
homogeneous saturated ideal defining a closed subscheme $Q_A \subseteq \PP^3_A$ that is
flat over $A$. By construction, $(Q_A) \times_{\Spec(A)} \Spec (\kappa (e))$ is the
unique quadric containing $X_e$, for all $e \in U$. Now, by the universal property of the
Hilbert scheme, there is a canonical morphism $f_U : U \to G$ and moreover, letting $U$
vary in an open affine covering of $E$, the morphisms $f_U$ glue together and produce the
required morphism $f$.

(b) Since double planes correspond to a closed subset of
$G$ the conclusion follows immediately from (a).
\end{proof}

The next Lemma shows the existence of the family we are interested in and
states some of
its  elementary properties.

\begin {lemma}\label{families_lemma2} Assume $d\ge 7$ and $g \le
\binom{d-3}{2}+1$. Then
we have:
\begin{itemize}
\item[(a)] There exists a reduced locally closed subscheme $\mathcal
F_{SET} \subset H_{d,g}$ which parameterizes the curves of subextremal type.
\item[(b)] There are closed subschemes $\mathcal F_{SE}$ and  $\mathcal
F^{(2)}_{SET}$ of $\mathcal F_{SET}$ which parameterize the subextremal
curves and the
curves of subextremal type lying in some double plane, respectively.
Moreover $\mathcal
F_{SET} = \mathcal F^{(2)}_{SET} \cup \mathcal F_{SE}$.
\end{itemize}
\end{lemma}

\begin{proof} By the definition of curve of subextremal type and by
semicontinuity it
follows that the subset corresponding to the curves of subextremal type is
locally
closed, hence it carries a natural structure of reduced subscheme of
$H_{d,g}$. This is
$\mathcal F_{SET}$.

(b) By Theorem \ref{rao_functions_subextremal_type}, a curve $C$ of
subextremal type is
subextremal if and only if $\rho_C(d+r-1) > 0$. Hence by semicontinuity
$\mathcal F_{SE}$
is a closed subset of $\mathcal F_{SET}$. Moreover the subset $\mathcal
F^{(2)}_{SET}$ is
closed by Lemma \ref{families_lemma1}. The last assertion is clear.
\end{proof}

Now we are going to study the two subfamilies $\mathcal F_{SE}$ and
$\cF^{(2)}_{SET}$ in
more detail. In particular we will show that, if $r\ge 3$, they are the two
irreducible
components of $\mathcal F_{SET}$ and have the same dimension.

We begin with $\cF^{(2)}_{SET}$. For this we will need the following concept:

\begin{definition} Let $b$ be an integer satisfying
$0 \le b \le \lfloor \frac{r-1}{2}\rfloor$. A curve $C$ is said to be {\it
of subextremal
type $b$} if $\rho_C = \rho_b$ (see
Theorem \ref{rao_functions_subextremal_type}). Thus a subextremal curve is
a curve of
subextremal type $0$.
\end{definition}
\medskip

\begin {theorem}\label{families1} Assume $d\ge 7$ and $r \ge 0$. Then:
\begin{itemize}
\item[(a)] The subscheme $\mathcal F^{(2)}_{SET} \subset H_{d, g}$ is irreducible and has
dimension $2r+6+\frac{(d-2)(d+1)}{2}$.
\item[(b)] If $r \geq 3$ there exists a stratification of $\mathcal
F^{(2)}_{SET}$ given by the Rao function, namely the curves of subextremal
type $b$ form
a non-empty irreducible subfamily ${\mathcal F}^{(2)}_{{SET}_b}$, which is:
\begin{itemize}
\item[(i)] open, generically smooth for even $r$ and smooth for odd $r$,
if $b=\lfloor \frac{r-1}{2}\rfloor$,
\item[(ii)] smooth, locally closed of codimension $1$, if
$0 < b < \lfloor \frac{r-1}{2}\rfloor$,
\item[(iii)] smooth, closed  of codimension $1$, if $b=0$
(it parameterizes the subextremal curves that are contained in some double
plane);
\end{itemize}
\end{itemize}
\end{theorem}

\begin{proof} We use the techniques of \cite {HS}. Let $2H$ be a fixed
double plane
and let $C\subseteq 2H$ be a curve of subextremal type. Then one can
associate to $C$ a flag of subschemes of $H$, namely $Z \subseteq C'
\subseteq D$, where
$D$ is the planar subcurve of $C$ of degree
$d-2$ (Theorem \ref {characterize_subextremal_type},(iv)) and $Z$ and $C'$
come from the
exact residual sequence with respect to $H$ (see Corollary \ref
{properties_subextremal_type},(ii)). Recall that $\deg Z = r$ and $\deg C'
= 2$.
Consider now the set of all curves of subextremal type contained in
$2H$ of degree $d$ and genus $g$. This coincides set theoretically with the
scheme
$H_{r,2,d-2}(2H)$ defined in \cite{HS}, which is irreducible and
generically smooth of
dimension $2r+3+\frac{(d-2)(d+1)}{2}$ (see \cite {HS}, Corollary 4.3).

Now by Lemma \ref{families_lemma1} there is a surjective morphism $f :
\cF^{(2)}_{SET} \to T$, where $T \cong \PP^3$ parameterizes the double planes.
Clearly the fibers of $f$ are homeomorphic to $H_{r,2,d-2}(2H)$. Then (a)
follows.

Now we prove (b). First of all observe that, due to the particular shape of
the Rao
functions (Theorem \ref{rao_functions_subextremal_type}) and by
semicontinuity, the
subsets $\mathcal F^{(2)}_{{SET}_b}$ are locally closed (open for $b = \lfloor
\frac{r-1}{2}\rfloor$, closed for $b=0$) and form a stratification of $\mathcal
F^{(2)}_{SET}$.

To prove the remaining properties it is easy to see that, by using the
morphism $f$
defined above, it is sufficient to study the problem in a fixed double
plane $2H$.

 According to \cite {HS} (see Propositions 4.1 and 4.2 and their proofs)
there is a smooth fibration $\pi: H_{r,2,d-2}(2H) \to D_{r,2}(H)$, where
$D_{r,2}(H)$ is
the flag scheme consisting of the pairs $(Z,C')$ where $C' \subseteq H$ is
a conic and
$Z$ is a locally complete intersection zero-dimensional scheme of degree
$r$ contained in
$C'$.

For each integer $b$ with $0 \le b \le \lfloor \frac{r-1}{2}\rfloor$, set
  $T_b(H) = \{(Z,C') \in D_{r,2}(H) \mid \partial h_Z = h_b\}$ (see
\S \ref{prel}) and let ${\mathcal F^{(2)}_{SET}}_b(H)$ be the set of curves of
subextremal type $b$ contained in $2H$. Then ${\mathcal F^{(2)}_{SET}}_b(H) =
\pi^{-1}(T_b(H))$ (see proof of Corollary
\ref{properties_subextremal_type}). Recall that
$D_{r,2}$ is irreducible and generically smooth of dimension $r+5$ (this
follows from
\cite {BH}: see \cite {HS}, proof of Proposition 4.2, for the idea of the
proof).
Moreover by \cite {F}, Theorem 2.4, the stratification of $D_{r,2}$ given
by the Betti
numbers has locally closed smooth irreducible strata. Now by Proposition
\ref{min_free_of_set} the Betti numbers of $Z$ depend only on $b$  if $0\le b <
\frac{r}{2}$, whence ${\mathcal F}^{(2)}_{{SET}_b}(H)$ is smooth and
irreducible for such
values of $b$. If $r$ is even and $b = \frac{r}{2}-1$ we have two different
strata, so we
can only say that ${\mathcal F}^{(2)}_{{SET}_ {\frac{r}{2}-1}}(H)$ is
irreducible and
generically smooth.

This proves (i) and the smoothness statements in (ii) and (iii).

Now we compute the dimensions of the strata. As above it is sufficient to
study the
problem for a fixed plane $H$.

If $b = 0$, then  $Z$ is a complete intersection $(1,r)$. These complete
intersections
form a smooth irreducible family of dimension $2+r$, as it is easily seen,
and hence
$\dim T_0 = 4+r = \dim D_{r,2}(H) - 1$. It follows that $\mathcal
F^{(2)}_{{SET}_0}(H)$
is smooth and $\dim {\mathcal F}^{(2)}_{{SET}_0}(H) = \dim H_{r,2,d-2}(2H)
- 1 = \dim
{\mathcal F}^{(2)}_{SET} -1$. This completes the proof of (iii).

It remains to compute $\dim ({\mathcal F}^{(2)}_{{SET}_b})$ for $0 < b <
\lfloor\frac
{r-1}{2}\rfloor$. For this it is sufficient to compute $\dim (T_b(H))$. Now
for $b$ in
the given range we have $r \ge 5$ and $r-b \ge 3$. Hence $\partial h_Z =
h_b$ if and only
if there is a line $L\subseteq H$ such that $\deg (L\cap Z) = r-b$ whence
the unique
conic containing $Z$ is reducible.

 This implies that
$\dim(T_b(H)) \le \dim (D_{r,2}(H)) - 1 = r+4$. Now the reduced schemes $Z
\subseteq H$
consisting of $r-b$ points on a line and $b$ points on a different line
form a non-empty
open subset $U \subseteq T_b(H)$. It is easy to see that $\dim (U) = r+4$
whence $\dim
(T_b(H)) = r+4$. It follows that $\dim {\mathcal F}^{(2)}_{{SET}_b} = \dim
{\mathcal
F}^{(2)}_{SET}-1$ for $b$ in the given range. This completes the proof of (ii).
\end{proof}

Now we turn our attention to the family of subextremal curves $\mathcal
F_{SE}$.

\begin {theorem}\label{families2} If $d\ge 7$  the subscheme
$\mathcal F_{SE} \subset H_{d,g}$ is irreducible. If moreover $r \ge 3$, then  $\mathcal
F_{SE}$ is smooth of  dimension $2r+6+\frac{(d-2)(d+1)}{2}$
\end{theorem}

\begin{proof} Let $C$ be a subextremal curve. By definition,
its Rao function $\rho_C$ depends only on $d$ and $g$ (see \S \ref{prel}). Moreover, by
Corollary \ref{Hilbert} we have that $\gamma_C$ can be computed in terms of $\rho_C$,
hence it is independent of the particular curve $C$. It follows that  the subextremal
curves are parameterized by the subscheme $H_{\gamma, \rho} \subseteq H_{d,g}$, as
defined in \cite {MDP1}, Definition VI.3.14, where $\rho:= \rho_C$ and $\gamma :=
\gamma_C$. Then we have $\mathcal F_{SE} = (H_{\gamma,\rho})_{red}$.

Now we show that $H_{\gamma, \rho}$ is irreducible. To this end we use some
ideas from liaison theory. Let
$$
0 \to F  \to N \oplus G \to I_C \to 0
$$
be the minimal N-type resolution of the subextremal curve $C$ where $F, G$
are free
$R$-modules of smallest possible rank. Then $N$ is the second syzygy
module of the Hartshorne-Rao module $M_C$ of $C$ (see \cite{MDP1}). Since
$M_C$ is a
Koszul module,  we know the minimal free resolution of $N$. Using the
mapping cone
procedure, the above sequence provides a free resolution of $I_C$.
Comparing with the
graded Betti numbers of $C$ (cf.\ Proposition \ref{prop-res-se}), we see
that we must
have $G = R(-2)$ and $F = R(r-4)\oplus R(-3) \oplus R(-d+1)$.
Hence, the corresponding modules in the minimal N-type resolutions of each two
subextremal curves $C$, $\widetilde{C}$ are isomorphic. Thus we conclude as
in Step (IV)
of the proof of \cite{N-gorliaison}, Theorem 7.3, that $C$ and
$\widetilde{C}$ belong to a flat family whose members belong to $H_{\gamma,
\rho}$ and
that is parameterized by an open subset of $\AAA^1$. The irreducibility of
$H_{\gamma,\rho}$ follows.

Now assume $r\ge 3$. By Theorem \ref{families1}(iii), the subextremal curves contained in
some double plane form an irreducible family of dimension $2r + 5 + \frac
{(d-2)(d+1)}{2}$. Since there are subextremal curves that are not contained in a double
plane (for example, perform a basic double linkage on a reduced reducible quadric
starting from an extremal curve of degree $d-2$ and genus $g-d+3$), we have
$$\dim \mathcal F_{SE} \geq 2r + 6 + \frac {(d-2)(d+1)}{2}.$$

Hence, to conclude our proof it is sufficient to show that the tangent space of
$H_{\gamma, \rho}$ at every closed point $t$ has dimension $t_{\gamma.\rho} = 2r + 6 +
\frac {(d-2)(d+1)}{2}$.

Let  $C$ be the curve corresponding to $t$ and let $M$ be the Rao module of
$C$. Then by
\cite {MDP1}, Theorem IX.4.2, the dimension of the tangent space of
$H_{\gamma, \rho}$ at
$t$ is
$$
t_{\gamma,\rho} = \delta_\gamma + \epsilon_{\gamma,\rho} - \dim_k
(\Hom(M,M)^0) +
\dim_k(\Ext^1(M,M)^0),
$$
where $\delta_\gamma$ and $\epsilon_{\gamma,\rho}$ are the number defined
in \cite
{MDP1}, ch. IX, 3.1.

The calculations are lengthy but elementary and make use of the assumption
$r\ge 3$. First of all one computes $\gamma$ from Corollary \ref {Hilbert},
and from this
one gets
$$\delta_\gamma = \frac {(d-2)(d+1)}{2} + 8 - r.$$

Next, from $\rho$ and $\gamma$ one finds $\epsilon_{\gamma,\rho} = r - 4$.
Since $M$ is a
Koszul module, one has that $\dim_k (\Hom(M,M)^0) = 1$ and from $\rho$ and
\cite {MDP1},
ch. IX, example 6.1, one gets $\dim_k(\Ext^1(M,M)^0) = 2r + 3$. It follows
$$t_{\gamma,\rho} = 2r + 6 + \frac {(d-2)(d+1)}{2}.
$$
This completes the proof.
\end{proof}

From Lemma \ref{families_lemma2} and Theorems \ref{families1} and \ref{families2} we have
immediately:

\begin {corollary}\label{families3} Assume $d \ge 7$ and
$r\geq 3$. The reduced subscheme $\mathcal F_{SET} \subset H_{d, g}$ is of
pure dimension
$2r+6+\frac{(d-2)(d+1)}{2}$ and its irreducible components are $\mathcal
F^{(2)}_{SET}$
and $\mathcal F_{SE}$. Moreover $\left (\mathcal F^{(2)}_{SET}\cap \mathcal
F_{SE}\right
)_{red}=\mathcal F^{(2)}_{{SET}_0}$.
\end{corollary}

\begin{remark} In the previous results we have made the assumption $r \ge 3$.
If $1\le r \le 2$ then $\cF_{SET} = \cF_{SE}$ by Theorem
\ref{rao_functions_subextremal_type} and the stratification of Theorem \ref{families1}(b)
is trivial. By Theorem \ref{families1}(a) and Theorem \ref {families2} we have that
$\cF_{SET}$ is irreducible and has dimension $\dim (\cF_{SET}) >
2r+6+\frac{(d-2)(d+1)}{2}$. The interested reader might carry out the calculation of
$t_{\gamma,\rho}$ as in the proof of Theorem \ref {families2} and get some more precise
information.
\end{remark}
\bigskip

\section{Two components of the Hilbert scheme} \label{sec-hilb}

In this section we show that, if $d\ge 7$ and $g < 0$, then the closures of
$\cF_{SE}$
and of $\cF_{SET}^{(2)}$ in $H_{d, g}$ are, topologically, irreducible
components of
$H_{d ,g}$. We use the same notation as in the previous section.

We begin with a geometrical description of the ``general subextremal curve.''

\begin{theorem} \label{general_subextremal_curves} Assume $d\ge 7$ and
$r\ge 3$. Then:
\begin{itemize}
\item[(a)] If $C'$ is an extremal curve of degree $d-1$ and genus $g-1$ and
$L$ is a $2$-secant line of $C'$ then $C : = L \cup C'$ is a sub-extremal
curve of degree
$d$ and genus $g$.
\item[(b)] There is a non-empty open set $U \subseteq \cF_{SE}$ such that
every curve $C \in U$ is as in {\rm (a)}.
\item[(c)] There is a non-empty open set $U' \subseteq U$ such that every $C
\in U'$ has a scheme-theoretical decomposition
$$C = D \cup Y \cup L$$
where $D$ is a smooth planar curve of degree $d-3$ spanning a plane $H$, $Y
\not
\subseteq H$ is a double line whose support lies in $H$ and $L \not
\subseteq H$ is a
$2$-secant line to $Y$. Moreover, the arithmetic genus of $Y$ satisfies
$g_Y \le - r$ and
$C_{red}$ is a curve of degree $d-1$ of maximal genus.
\end{itemize}
\end {theorem}

\begin{proof} (a) From \S 2.2 it follows that the maximum of
$\rho_{C'}$ is $r$ and $C'$ contains a subcurve $P$ of degree $d-2$
spanning a plane $H$
and that the residual curve of $C'$ with respect to $H$ is a line $\ell$.
It is clear
that $\deg C = d$ and by the genus formula for the union of two curves it
follows that
$p_a(C) = g$. Since $\deg P \ge 3$ we have $L \not \subseteq H$, whence $C$
contains a
planar subcurve of degree $d-2$, namely $P$, but it does not contain a
planar subcurve of
degree $d-1$.

 Now by \cite{MDP3}, Proposition 0.6 and our numerical
assumptions, we have that $\ell \subseteq P$ and $C = D'\cup Q$, where $Q$
is a multiple
line supported by $\ell$ and $D' \not \supseteq \ell$. Then it is clear
that $L$ must be
a secant line of $Q$ and, in particular, $L$ meets $\ell$. This implies
that the residual
curve of $C$ with respect to $H$ is the planar degree $2$ curve $\ell \cup
L$. Hence $C$
is of subextremal type by Theorem \ref{characterize_subextremal_type}.
Moreover since
$L\not \subseteq H$ the unique quadric containing $C$ is reduced whence $C$ is
subextremal by Theorem \ref{rao_functions_subextremal_type}. This proves (a)

(b) Now we want to show that the subextremal curves constructed above form
a family,
and we want to compute its dimension.

Let $E \times  \PP^3\supseteq X\to E$ be the family of extremal curves of
degree $d-1$
and genus $g-1$ and let $G \times \PP^3\supseteq Y\to G$ be the
Grassmannian of lines of
$\PP^3$. Recall that $\dim E = 2r + 4 + \frac {(d-2)(d+1)}{2}$ by Theorem
\ref{family extremal curves}.

Now the family
$$
E \times G \times \PP^3\supseteq (X \times G) \cap
 (E \times Y) \to E \times G
$$
parameterizes bijectively the intersections $X_e\cap Y_g$ and by
Chevalley's theorems
there is a locally closed subset $V\subseteq E \times G$ such that $(e,g)\in V$
if and only if $\length (X_e\cap Y_g) = 2$  (that is if and only if $Y_g$
is a $2$-secant
line of $X_e$).

For any $e \in E$ let $H_e$ be the plane containing the planar subcurve of
$X_e$ of
degree $d-2$ and let $\ell_e$ be the residual line of $X_e$ with respect to
$H_e$. Let
$Q_e$ be the largest subcurve of $X_e$ supported by $\ell_e$. Then, as we
have seen
above, $Y_g$ is a $2$-secant line of $X_e$ if and only if it is a
$2$-secant line of
$Q_e$. Now it easily follows that the fibers of the projection $V\to E$
have dimension
$2$, whence $\dim V = \dim E + 2 = 2r + 6 + \frac{(d-2)(d+1)}{2}$.

Consider now the family
$$ E \times G \times \PP^3\supseteq (E \times Y)  \cup
 (X \times G) \overset \varphi \to E \times G.$$

It is easy to see that it parameterizes bijectively the schemes $X_e\cup
Y_g$. Thus, it
follows that the subextremal curves constructed in (a) are exactly the
curves of the
family $\varphi^{-1}(V) \to V$.

By the universal property of the Hilbert scheme there is an injective morphism
$\Phi : V \to H_{d,g}$ and $\Phi(V) \subseteq \mathcal F_{SE}$ by (a). Moreover
$\Phi(V)$ is constructible and since $\dim V = \dim  \mathcal F_{SE}$ and
$\mathcal
F_{SE}$ is irreducible by Theorem \ref{families2} it follows that $\Phi(V)$
contains a
non-empty open subset $U$ of $\mathcal F_{SE}$ and the conclusion follows.

(c) By \cite{MDP3}, Proposition 0.6 it follows that there is a non-empty open subset $E'
\subseteq E$ such that every $C' \in E'$ has a scheme-theoretical decomposition $C' = D
\cup Y$, where $D$ is a planar smooth curve of degree $d-3$ spanning a plane $H$ and $Y$
is a double line whose support lies in $H$. Let $\pi : E\times G \to E$ be the
projection. Then one shows as above that the image of $\pi^{-1}(E')\cap V$ in $H_{d,g}$
contains a non-empty open subset $U'$ of $\cF_{SE}$ with the required properties.

The genus of $Y$ can be easily bounded by using the formula $g = p_a(D) +
p_a(Y) +
\length D \cap Y - 1$, observing that $\length (D \cap Y) \ge d-3$.

Finally if $C \in U'$ then $C_{red} = D \cup Y_{red}\cup L$ is a
curve of degree $d-1$ and of maximal genus, being the union of a planar
curve and a
line meeting it in a scheme of length $1$ (see \cite {H}).
\end {proof}

Now we can deal with $\cF_{SE}$.

\begin{theorem}\label{subextremal_are_component} Assume $d \ge 7$ and $g <
0$. Then the
closure of $\cF_{SE}$ in $H_{d,g}$ is, topologically, an irreducible
component of
$H_{d,g}$.
\end{theorem}

\begin {proof} Observe first that our numerical assumptions imply  $r \ge
8$, whence, in
particular, Theorem \ref{general_subextremal_curves} applies.

We follow some ideas from \cite{MDP3}, proof of Proposition 3.6. By Theorem
\ref{families2}, there is   an irreducible component $\cF$ of $H_{d,g}$
containing
$\cF_{SE}$ We want to show that $\cF = \overline{\cF_{SE}}$. We argue by
contradiction,
assuming that $\cU :=\cF\setminus \overline{\cF_{SE}} \ne \emptyset$.

We use the scheme-theoretical decomposition $C_0 = D \cup Y \cup L$ of a
general $C_0 \in
\cF_{SE}$  given in  Theorem \ref {general_subextremal_curves} (c).

Fix $C_0$ as above and let $C\in \cU$. Then there is a flat family $X \to
T$, where $T$
is a smooth connected curve, and points $t_0, u \in T$ such that $C_0 =
X_{t_0}$,
$C = X_u$, and $X_t \in \cU$ for $t\in T\setminus \{t_0\}$.

By our assumption on $g$ every curve $C \in \cU$ is non-integral, hence it is
either non-reduced or reduced and reducible. We consider the two cases
separately.

\underline {Case 1}. Assume that $C$ is non-reduced. Then $X$ is non
reduced and we set
$X':= X_{red}$. Since $T$ is a smooth curve the family $X'\to T$ is flat.
Moreover $X'_t
= (X_t)_{red} \neq X_t$ for general $t\in T$ and $(C_0)_{red} \subseteq
X'_{t_0}
\subsetneq C_0$. By the particular shape of $C_0$ we have  $\deg
(C_0)_{red} = \deg C_0
-1 = d-1$, whence $(C_0)_{red} = X'_{t_0}$. By flatness we have that $X'_t$
is a curve of
degree $d-1$ and maximal genus, as $(C_0)_{red}$ is such a curve. It
follows that $X'_t$
is the union of a planar curve $P_t$ and of a line $\ell_t$ meeting $P_t$
in a scheme of
length 1 (see \cite {H}). Now by degree reasons $X_t$ contains a double
line $Y_t$. If
$(Y_t)_{red} \subseteq P_t$ we have that $X_t\in \cF_{SE}$ whence $C \in
\cF_{SE}$, a
contradiction. So we have $(Y_t)_{red} = \ell_t$. We want to show that this
leads again
to a contradiction.

Assume first that $P_t$ is integral. Then $X'_t$ has two irreducible
components, namely
$P_t$ and $\ell_t$. It follows that $X'$ has two irreducible components
$X'_1$ and $X'_2$
corresponding to $P_t$ and $\ell_t$, respectively. Since $X$ has no
embedded components,
it follows that it has exactly two irreducible components $X_1$ and $X_2$
that are, by
degree reasons, topologically equal to $X'_1$ and $X'_2$, respectively.
This implies that
$(X_2)_{t_0}$ is a sub-curve of $C_0$ supported by $L$, hence $(X_2)_{t_0}
= L$. But this
is a contradiction because the families $X_i \to T$ are flat, being $T$ a
smooth curve,
while $\deg (X_2)_{t_0} \ne \deg (X_2)_t$. Hence, $P_t$ is not integral.

It follows that every general $C \in \cU$ has a scheme-theoretical
decomposition $C = P
\cup W$, where $P$ is a non-integral planar curve of degree $d-2$ and $W$
is a double
line whose support meets $P$ but does not lie in the plane spanned by $P$.
In particular
$\epsilon:= \length W \cap P \in \{1,2\}$.

Now since $g = p_a(P) + p_a(W) + \epsilon -1$ we get $-r+1 \ge p_a(W) \ge
-r$, whence, in
particular, $p_a(W) \le -2$. It also follows that the double lines $W$ move
in a family
of dimension $5-2p_a(W) \le 5+2r$ (see \cite {MDP3}, Theorem 4.1). Since
the non-integral
planar curves of degree $d-2$ move in a family of dimension $\frac 12
d(d-3) + 5$ we have
$\dim \cF = \dim \cU \le 5+2r + \frac 12 d(d-3) + 5 \le 2r + 6 +\frac 12
(d-2)(d+1) =
\dim \cF_{SE}$ where the last equality is due to  Theorem \ref{families2}.
This is a
contradiction, and the conclusion follows.

\underline {Case 2}. Assume that $C$ is reduced. Then $X$ is reduced and
reducible, namely there is a proper scheme decomposition $X = X_1 \cup X_2$
and the
families $X_i \to T$ are flat since $T$ is a smooth curve. Moreover we
have, set-theoretically, $(X_1)_{t_0}\cup (X_2)_{t_0} = C_0$.
Up to interchanging $X_1$ and $X_2$ we have three possibilities, namely:
\begin{itemize}
\item[(i)] $((X_1)_{t_0})_{red} = L$ and $((X_2)_{t_0})_{red} = Y_{red}
\cup D$
\item[(ii)] $((X_1)_{t_0})_{red} = Y_{red}$ and $((X_2)_{t_0})_{red} = L
\cup D$
\item[(iii)] $((X_1)_{t_0})_{red} = D$ and $((X_2)_{t_0})_{red} = Y_{red}
\cup L$
\end{itemize}

If (i) holds then $(X_1)_{t_0} = L$ whence, by degree reasons, $(X_2)_{t_0} =
Y \cup D$. It follows that $(X_2)_{t_0}$ is an extremal curve of degree
$d-1$ and
genus $g-1 < 0$. Then $(X_2)_t$ is non-integral, hence it is extremal by
\cite {MDP3},
Proposition 3.6. Now (as in the proof of Theorem
\ref{general_subextremal_curves}, (a))
$r$ is the maximum of the Rao function of $(X_2)_{t_0}$ and hence of
$(X_2)_t$. Since
 $r \ge 8$, $(X_2)_t$ is not reduced (see \cite{MDP3}, Proposition 0.6),
 whence $X$ is not reduced, a contradiction.

If (ii) holds then $(X_2)_{t_0} = L \cup D$, whence $(X_1)_{t_0} = Y$. Then
$(X_1)_t$ is
a curve of degree $2$. Moreover by Theorem
\ref{general_subextremal_curves},(c) we have
$p_a((X_1)_t) = p_a(Y) \le -r \le -3$, which implies that $(X_1)_t$ is not
reduced,
whence $X$ is not reduced, again a contradiction.

If (iii) holds we get, arguing as above, $(X_2)_{t_0} = L \cup Y$. It is
easy to show
that $p_a(L \cup Y) \le -r+1 \le - 7$. Since $\deg (L \cup Y) = 3$ we get
$p_a((X_2)_t)
\le - \deg((X_2)_t)$, whence $X_t$ is not reduced. As above, it follows
that $X$ is not
reduced, a contradiction.
\end{proof}

Now we consider $\cF_{SET}^{(2)}$. Our strategy is similar to the previous
one for
$\cF_{SE}$. We begin with a geometric description of the general curve in
$\cF_{SET}^{(2)}$.

\begin{lemma} \label{general_set_curves} Assume $d\ge 7$. Then there is a
non-empty open set  $U \subseteq \cF_{SET}^{(2)}$ such that every $C \in U$
admits a
scheme-theoretical decomposition $C = Y \cup E$, where $E$ is a smooth
planar curve of
degree $d-4$ contained in a plane $H$ and $Y$ is a curve of degree $4$
whose support is a
smooth conic contained in $H$. Moreover, $p_a(Y) \le -r+1$.
\end{lemma}

\begin{proof} Let $H$ be a fixed plane. Then by \cite{HS} there is a morphism
$\sigma: H_{r,2,d-2}(2H) \to D_{r,2,d-2}(H)$, where $D_{r,2,d-2}(H)$ is the
flag scheme
parameterizing the triples $(Z,C',D)$ (see proof of Theorem
\ref{families1}). The subset
 $D'_{r,2,d-2}(H)\subseteq D_{r,2,d-2}(H)$,  where
$C'$ is a smooth conic, $D = C' \cup E$, and $E$ is smooth, is open and
non-empty, thus
$\pi^{-1}(D'_{r,2,d-2}(H))$ is a non-empty open subset of $H_{r,2,d-2}(2H)$
whose curves
have the required scheme-theoretical decomposition.

Let now $P\in \PP^3$ be a point and denote by $(\cF_{SET}^{(2)})_P$ the set
of curves in
$\cF_{SET}^{(2)}$ lying in a double plane not containing $P$. By Lemma
\ref{families_lemma1} the subsets $(\cF_{SET}^{(2)})_P$, as $P$ varies,
form an open
covering of $\cF_{SET}^{(2)}$. For every $P\in \PP^3$ fix a plane $H_P$ not
containing
$P$. Then the projection from $P$ induces a surjective morphism $g_P :
(\cF_{SET}^{(2)})_P \to H_{r,2,d-2}(2H_P)_{red}$ and the existence of $U$
follows.

Let now $C = Y\cup E$ be a curve in $U$. Then the scheme $E\cap Y$ is
zero-dimensional
and $\length (E\cap Y) \ge 2(d-4)$. Since $g = p_a(Y) + p_a(E) + \length
(E\cap Y) - 1$,
the bound for $p_a(Y)$ is obtained by a straightforward calculation.
\end{proof}

\begin{theorem}\label{subextremal_type_2_are_component} Assume $d \ge 7$
and $g < 0$.
Then the closure of $\cF_{SET}^{(2)}$ in $H_{d,g}$ is, topologically, an
irreducible
component of $H_{d,g}$.
\end{theorem}

\begin{proof} We use the same setting and notation (with
obvious modifications) as in the proof of Theorem \ref
{subextremal_are_component}. In
particular $C_0 = E \cup Y$ will have the structure given by Lemma
\ref{general_set_curves}. Observe also that our numerical assumptions imply
$r \ge 8$,
whence $p_a(Y) \le -7$.

\underline {Case 1}. The general curve in $\cU$ is not reduced. Then, $Y$ being
irreducible, we have $X'_{t_0} = E \cup Y_{red}$. Thus there are only two
possibilities,
namely:
\begin{itemize}
\item[(i)] $X_t$ contains a double conic; or
\item[(ii)] $X_t$ contains a multiple line.
\end{itemize}
If (i) holds then $X_t \in \cF_{SET}^{(2)}$ and we are done. If (ii) holds
then $X'_t$ is
the union of a curve of degree $d-3$ and a line. Let $X' = X'_1 \cup X'_2$
be the
corresponding decomposition. Then $(X'_2)_{t_0}$ is  a line contained in $E
\cup
Y_{red}$, which is impossible.

\underline {Case 2}. The general curve in $\cU$ is reduced and reducible.
then we have $X
= X_1 \cup X_2$ with $(X_1)_{t_0} = Y$. But $p_a((X_1)_t) =
p_a((X_1)_{t_0}) \le -7 < -4
= -\deg ((X_1)_t)$. Then $X$ is not reduced, a contradiction.
\end{proof}

With respect to the Hilbert scheme, our results about curves of subextremal
type can be
summarized as follows:

\begin{corollary} \label{cor-hilb}
If $d \geq 7$ and $g < 0$, then the Hilbert scheme $H_{d, g}$ has at least
three
components. Topologically, two components are the closures of
$\cF^{(2)}_{SET}$ and
$\cF_{SE}$ and another component is formed by the closure of $\cF_{EX}$ that
parameterizes the extremal curves. The first two components have dimension
$\frac{3}{2} d
(d-5) + 19 - 2g$ and meet in a subscheme of codimension one, the third
component has
dimension
$\frac{3}{2} d (d - 3) + 9 - 2g$. The support of these three components is
generically
smooth.
\end{corollary}

\begin{proof} It suffices to note that the results about $\cF_{EX}$ are
shown in \cite{MDP3}.
\end{proof}

We believe that the above result remains true if we replace the assumption
on the genus
by $g \leq \binom{d-3}{2}-2$. However, proving the statement in this
generality  seems to
require a different approach.

\end{document}